\documentclass[reqno,11pt]{amsart}
\usepackage{amsfonts}

\usepackage{graphicx}
\usepackage{pstricks}
\usepackage{amsmath}
\usepackage{amsmath}
\usepackage{amsxtra}

\theoremstyle{plain}
\newtheorem{theorem}{Theorem}[section]
\newtheorem{lemma}[theorem]{Lemma}

\newtheorem{proposition}[theorem]{Proposition}

\newtheorem{conjecture}[theorem]{Conjecture}
\newtheorem{definition}[theorem]{Definition}

\newtheorem{problem}[theorem]{Problem}

\theoremstyle{definition}

\newtheorem{remark}[theorem]{Remark}

\numberwithin{equation}{section} 
\begin{document}
\title{Reconstruction of the Berger measure when the core is of tensor form}
\author{Ra\'{u}l E. Curto}
\address{Department of Mathematics, The University of Iowa, Iowa City, Iowa
52242, USA}
\email{rcurto@math.uiowa.edu}
\urladdr{http://www.math.uiowa.edu/\symbol{126}rcurto/}
\author{Sang Hoon Lee}
\address{Department of Mathematics, Chungnam National University, Daejeon,
305-764, Korea}
\email{shlee@math.cnu.ac.kr}
\urladdr{}
\author{Jasang Yoon}
\address{Department of Mathematics, Iowa State University, Ames, Iowa 50011,
USA}
\email{jyoon@iastate.edu}
\urladdr{http://www.public.iastate.edu/\symbol{126}jyoon/}
\thanks{The first named author was partially supported by NSF Grants
DMS-0099357 and DMS-0400471.}
\subjclass[2000]{Primary 47B20, 47B37, 47A13, 28A50; Secondary 44A60, 47-04,
47A20}
\keywords{jointly hyponormal pairs, subnormal pairs, $2$-variable weighted
shifts, Berger measure, propagation phenomena, flatness}

\begin{abstract}
Let $\mathfrak{H}_{0}$ denote the class of commuting pairs of subnormal
operators on Hilbert space, and let $\mathcal{TC}:=\{\mathbf{T}\in \mathfrak{%
H}_{0}:c(\mathbf{T)}$ is of tensor form$\}$, where $c(\mathbf{T})$ is the
core of $\mathbf{T}$. \ We obtain a concrete necessary and sufficient
condition for the subnormality of $\mathbf{T}\equiv (T_{1},T_{2})\in 
\mathcal{TC}$ in terms of $c(\mathbf{T})$, the marginal measures of $T_{1}$
and $T_{2}$, and the weight $\alpha _{01}$.
\end{abstract}

\maketitle

\section{\label{Int}Introduction}

The Lifting Problem for Commuting Subnormals (LPCS) asks for necessary and
sufficient conditions for a pair of subnormal operators on Hilbert space to
admit commuting normal extensions. \ It is well known that the commutativity
of the pair is necessary but not sufficient (\cite{Abr}, \cite{Lu1}, \cite%
{Lu2}, \cite{Lu3}), and it has recently been shown that the joint
hyponormality of the pair is necessary but not sufficient \cite{CuYo1}. \
Abstract solutions of LPCS were given in \cite[Theorem 3.1]{CLY1} and %
\cite[Theorem 2.7]{Yo1}, while concrete, necessary conditions, albeit not
sufficient, for the lifting were found in \cite[Theorem 3.3 ]{CuYo2} and %
\cite[Theorem 2.10 ]{Yo1}) in the case of $2$-variable weighted shifts. \ In
(\cite{CuYo1}, \cite{CuYo2}, \cite{propagation}, \cite{CLY1}, \cite{power}, %
\cite{Yo1} and \cite{Yo2}) we have shown that many of the nontrivial aspects
of LPCS are best detected within the class $\mathfrak{H}_{0}$ of commuting
pairs of subnormal operators; we thus focus our attention on this class. \
More generally, we will denote the class of subnormal pairs by $\mathfrak{H}%
_{\infty }$, and for each integer $k\geq 1$ the class of $k$-hyponormal
pairs in $\mathfrak{H}_{0}$ by $\mathfrak{H}_{k}$. \ Clearly, $\mathfrak{H}%
_{\infty }\subseteq \cdots \subseteq \mathfrak{H}_{k}\subseteq \cdots
\subseteq \mathfrak{H}_{2}\subseteq \mathfrak{H}_{1}\subseteq \mathfrak{H}%
_{0}$; the main results in \cite{CuYo1} and \cite{CLY1} show that these
inclusions are all proper. \ The constructions in \cite{CuYo1} and \cite%
{CuYo2} have shed light on structural and spectral properties of
multivariable weighted shifts, and have brought about some new phenomena in
joint spectral theory. \ More recently, we have made use of the tools and
techniques in those papers and in \cite{propagation} and \cite{CLY1} to
approach LPCS from a new angle: to what extent the subnormality of the
powers of a $2$-variable weighted shift can detect the subnormality of the
pair. \ In \cite{power} we discovered a large class of $2$-variable weighted
shifts $\mathbf{T}\equiv (T_{1},T_{2})$ for which the subnormality of $%
(T_{1}^{2},T_{2})$ and $(T_{1},T_{2}^{2})$ does imply the subnormality of $%
\mathbf{T}$. \ This is the class $\mathcal{TC}$ (see Definition \ref{tc}
below). \ 

In this paper we study the subnormality of $2$-variable weighted shifts $%
\mathbf{T}\in \mathcal{TC}$. \ Since a $2$-variable weighted shift is
subnormal if and only if its weight moments are the moments of a probability
measure, known as the Berger-Gellar-Wallen measure (or briefly Berger
measure), the search for necessary and sufficient conditions leads to the
following concrete problem.

\begin{problem}
\label{tensorcore}Let $\mathbf{T\in }\mathcal{TC}$ and assume $\mathbf{T}$
is hyponormal. \ Additionally, assume that $c(\mathbf{T)}$ is subnormal,
with Berger measure $\xi \times \eta $. \ Find necessary and sufficient
conditions on the rest of the weight data to guarantee the subnormality of $%
\mathbf{T}$.
\end{problem}

Problem \ref{tensorcore} is a special instance of the
Reconstruction-of-the--Measure Problem, which we now describe. \ Given $%
\mathbf{T\in }\mathfrak{H}_{0}$, the $j$-th row and the $i$-th column of the
weight diagram have their own Berger measures, $\xi _{j}$ and $\eta _{i}$,
respectively. \ Solving LPCS in this case amounts to finding a measure $\mu $
on $\mathbb{R}_{+}^{2}$ which interpolates $\{\xi _{j},\eta
_{i}\}_{i,j=0}^{\infty }$. \ Without loss of generality we can assume that $%
\xi _{j+1}\ll \xi _{j}$ and $\eta _{i+1}\ll \eta _{i}\;($all $i,j\geq 0)$ %
\cite[Theorem 3.3]{CuYo2}. \ Moreover, $\xi _{j}$ must equal $\mu _{j}^{X}$
(the marginal measure of $\mu _{j}$), where $d\mu _{j}(s,t):=\frac{1}{\gamma
_{0j}}t^{j}d\mu (s,t)$; in fact, $d\xi _{j}(s)=\{\frac{1}{\gamma _{0j}}%
\int_{Y}t^{j}\;d\Phi _{s}(t)\}\;d\mu ^{X}(s)$, where $d\mu (s,t)\equiv d\Phi
_{s}(t)\;d\mu ^{X}(s)$ is the disintegration of $\mu $ by vertical slices %
\cite[Theorem 3.1]{CuYo2}; and similarly for $\eta _{i}$. \ From this
perspective, LPCS consists of ``compatibly gluing together''\ the measures $%
\xi _{j}$ and $\eta _{i}$ on $\mathbf{R}_{+}$ to produce a measure $\mu $ on 
$\mathbf{R}_{+}^{2}$ which satisfies the required properties to be the
Berger measure of $\mathbf{T}$. \ We claim this can be done explicitly for $%
\mathbf{T}\in \mathcal{TC}$.\ \ Note that $\xi _{0}=\mu ^{X}$ and $\eta
_{0}=\mu ^{Y}$. \ 

Our main result is Theorem \ref{thm1}, which provides a complete solution to
Problem \ref{tensorcore}: $\mathbf{T}\equiv (T_{1},T_{2})\in \mathcal{TC}$
is subnormal if and only if measures $\psi $ and $\varphi $ given by (\ref%
{eqpsi}) and (\ref{eqphi}), respectively, are positive. \ As an application,
we give a concrete condition for the subnormality of flat $2$-variable
weighted shifts (Proposition \ref{ex1}).

Let $\mathcal{H}$ be a complex Hilbert space and let $\mathcal{B}(\mathcal{H}%
)$ denote the algebra of bounded linear operators on $\mathcal{H}$. \ We say
that $T\in \mathcal{B}(\mathcal{H})$ is \textit{normal} if $T^{\ast
}T=TT^{\ast },$ and \textit{subnormal} if $T=N|_{\mathcal{H}}$, where $N$ is
normal and $N(\mathcal{H})$ $\mathcal{\subseteq H}$. \ An operator $T$ such
that $T^{\ast }T\geq TT^{\ast }$ is said to be \textit{hyponormal}. For $%
S,T\in \mathcal{B}(\mathcal{H})$ let $[S,T]:=ST-TS$.\ We say that an $n$%
-tuple $\mathbf{T:}=(T_{1},\cdots ,T_{n})$ of operators on $\mathcal{H}$ is
(jointly) \textit{hyponormal} if the operator matrix 
\begin{equation*}
\lbrack \mathbf{T}^{\ast },\mathbf{T]:=}\left( 
\begin{array}{llll}
\lbrack T_{1}^{\ast },T_{1}] & [T_{2}^{\ast },T_{1}] & \cdots & [T_{n}^{\ast
},T_{1}] \\ 
\lbrack T_{1}^{\ast },T_{2}] & [T_{2}^{\ast },T_{2}] & \cdots & [T_{n}^{\ast
},T_{2}] \\ 
\text{ \thinspace \thinspace \quad }\vdots & \text{ \thinspace \thinspace
\quad }\vdots & \ddots & \text{ \thinspace \thinspace \quad }\vdots \\ 
\lbrack T_{1}^{\ast },T_{n}] & [T_{2}^{\ast },T_{n}] & \cdots & [T_{n}^{\ast
},T_{n}]%
\end{array}%
\right)
\end{equation*}%
is positive on the direct sum of $n$ copies of $\mathcal{H}$ (cf. \cite{Ath}%
, \cite{CMX}). \ The $n$-tuple $\mathbf{T}$ is said to be \textit{normal} if 
$\mathbf{T}$ is commuting and each $T_{i}$ is normal, and $\mathbf{T}$ is 
\textit{subnormal }if $\mathbf{T}$ is the restriction of a normal $n$-tuple
to a common invariant subspace. \ The Bram-Halmos criterion for subnormality
states that an operator $T\in \mathcal{B}(\mathcal{H})$ is subnormal if and
only if 
\begin{equation*}
\sum_{i,j}(T^{i}x_{j},T^{j}x_{i})\geq 0
\end{equation*}%
for all finite collections $x_{0},x_{1},\cdots ,x_{k}\in \mathcal{H}$ (\cite%
{Bra}, \cite{Con}). \ Using Choleski's algorithm for operator matrices, it
is easy to see this is equivalent to the $k$-tuple $(T,T^{2},\cdots ,T^{k})$
is hyponormal for all $k\geq 1$.

For $\alpha \equiv \{\alpha _{n}\}_{n=0}^{\infty }$ a bounded sequence of
positive real numbers (called \textit{weights}), let $W_{\alpha }:\ell ^{2}(%
\mathbb{Z}_{+})\rightarrow \ell ^{2}(\mathbb{Z}_{+})$ be the associated
unilateral weighted shift, defined by $W_{\alpha }e_{n}:=\alpha
_{n}e_{n+1}\;($all $n\geq 0)$, where $\{e_{n}\}_{n=0}^{\infty }$ is the
canonical orthonormal basis in $\ell ^{2}(\mathbb{Z}_{+}).$ \ For notational
convenience, we will write sometimes \textquotedblleft $shift(\alpha
_{0},\alpha _{1},\cdots )$" for $W_{\alpha }$. \ In particular, $%
U_{+}:=shift(1,1,\cdots )$ and $S_{a}:=shift(a,1,1,\cdots ).$ \ For a
weighted shift $W_{\alpha }$, \textit{the moments of $\alpha $} are given as 
\begin{equation*}
\gamma _{k}(W_{\alpha })\equiv \gamma _{k}(\alpha ):=%
\begin{cases}
1 & \text{if }k=0 \\ 
\alpha _{0}^{2}\cdots \alpha _{k-1}^{2} & \text{if }k>0.%
\end{cases}%
\end{equation*}%
It is easy to see that $W_{\alpha }$ is never normal, and that it is
hyponormal if and only if $\alpha _{0}\leq \alpha _{1}\leq \cdots $. \
Similarly, consider double-indexed positive bounded sequences $\alpha _{%
\mathbf{k}},\beta _{\mathbf{k}}\in \ell ^{\infty }(\mathbb{Z}_{+}^{2})$ , $%
\mathbf{k}\equiv (k_{1},k_{2})\in \mathbb{Z}_{+}^{2}:=\mathbb{Z}_{+}\times 
\mathbb{Z}_{+}$ and let $\ell ^{2}(\mathbb{Z}_{+}^{2})$ be the Hilbert space
of square-summable complex sequences indexed by $\mathbb{Z}_{+}^{2}$. \
(Recall that $\ell ^{2}(\mathbb{Z}_{+}^{2})$ is canonically isometrically
isomorphic to $\ell ^{2}(\mathbb{Z}_{+})\bigotimes \ell ^{2}(\mathbb{Z}_{+})$%
.) \ We define the $2$-variable weighted shift $\mathbf{T}:=(T_{1},T_{2})$
by 
\begin{equation*}
\begin{cases}
T_{1}e_{\mathbf{k}}:=\alpha _{\mathbf{k}}e_{\mathbf{k+}\varepsilon _{1}} \\ 
T_{2}e_{\mathbf{k}}:=\beta _{\mathbf{k}}e_{\mathbf{k+}\varepsilon _{2}},%
\end{cases}%
\end{equation*}%
where $\mathbf{\varepsilon }_{1}:=(1,0)$ and $\mathbf{\varepsilon }%
_{2}:=(0,1)$. Clearly, 
\begin{equation}
T_{1}T_{2}=T_{2}T_{1}\Longleftrightarrow \beta _{\mathbf{k+}\varepsilon
_{1}}\alpha _{\mathbf{k}}=\alpha _{\mathbf{k+}\varepsilon _{2}}\beta _{%
\mathbf{k}}\;(\text{all }\mathbf{k}).  \label{commuting}
\end{equation}%
In an entirely similar way one can define multivariable weighted shifts.

Trivially, a pair of unilateral weighted shifts $W_{\alpha }$ and $W_{\beta
} $ gives rise to a $2$-variable weighted shift $\mathbf{T}\equiv
(T_{1},T_{2}) $, if we let $\alpha _{(k_{1},k_{2})}:=\alpha _{k_{1}}$ and $%
\beta _{(k_{1},k_{2})}:=\beta _{k_{2}}\;$(all $k_{1},k_{2}\in \mathbb{Z}_{+}$%
). In this case, $\mathbf{T}$ is subnormal (resp. hyponormal) if and only if
so are $T_{1}$ and $T_{2}$; in fact, under the canonical identification of $%
\ell ^{2}(\mathbb{Z}_{+}^{2})$ and $\ell ^{2}(\mathbb{Z}_{+})\bigotimes \ell
^{2}(\mathbb{Z}_{+})$, $T_{1}\cong I\bigotimes W_{\alpha }$ and $T_{2}\cong
W_{\beta }\bigotimes I$, and $\mathbf{T}$ is also doubly commuting. For this
reason, we do not focus attention on shifts of this type, and use them only
when the above mentioned triviality is desirable or needed. \ Given $\mathbf{%
k}\in \mathbb{Z}_{+}^{2}$, the moment of $(\alpha ,\beta )$ of order $%
\mathbf{k}$ is 
\begin{equation*}
\gamma _{\mathbf{k}}(\mathbf{T})\equiv \gamma _{\mathbf{k}}(\alpha ,\beta ):=%
\begin{cases}
1 & \text{if }\mathbf{k}=0 \\ 
\alpha _{(0,0)}^{2}\cdots \alpha _{(k_{1}-1,0)}^{2} & \text{if }k_{1}\geq 1%
\text{ and }k_{2}=0 \\ 
\beta _{(0,0)}^{2}\cdots \beta _{(0,k_{2}-1)}^{2} & \text{if }k_{1}=0\text{
and }k_{2}\geq 1 \\ 
\alpha _{(0,0)}^{2}\cdots \alpha _{(k_{1}-1,0)}^{2}\beta
_{(k_{1},0)}^{2}\cdots \beta _{(k_{1},k_{2}-1)}^{2} & \text{if }k_{1}\geq 1%
\text{ and }k_{2}\geq 1.%
\end{cases}%
\end{equation*}

(We remark that, due to the commutativity condition (\ref{commuting}), $%
\gamma _{\mathbf{k}}$ can be computed using any nondecreasing path from $%
(0,0)$ to $(k_{1},k_{2})$.) \ We now recall a well known characterization of
subnormality for multivariable weighted shifts \cite{JeLu}, due to C. Berger
(cf. \cite[III.8.16]{Con}) and independently established by Gellar and
Wallen \cite{GeWa}) in the single variable case: $\ \mathbf{T\equiv (}%
T_{1},T_{2})$ admits a commuting normal extension if and only if there is a
probability measure $\mu $ (which we call the Berger measure of $\mathbf{T}$%
) defined on the $2$-dimensional rectangle $R=[0,a_{1}]\times \lbrack
0,a_{2}]$ (where $a_{i}:=\left\| T_{i}\right\| ^{2}$) such that $\gamma _{%
\mathbf{k}}=\int_{R}s^{k_{1}}t^{k_{2}}d\mu (s,t),$ for all $\mathbf{k}\in 
\mathbb{Z}_{+}^{2}$. \ In the single variable case, if $W_{\alpha }$ is
subnormal with Berger measure $\xi _{\alpha }$ and $h\geq 1$, and if we let $%
\mathcal{L}_{h}:=\bigvee \{e_{n}:n\geq h\}$ denote the invariant subspace
obtained by removing the first $h$ vectors in the canonical orthonormal
basis of $\ell ^{2}(\mathbb{Z}_{+})$, then the Berger measure of $W_{\alpha
}|_{\mathcal{L}_{h}}$ is $\frac{s^{h}}{\gamma _{h}}d\xi (s)$; alternatively,
if $S:\ell ^{\infty }(\mathbb{Z}_{+})\rightarrow \ell ^{\infty }(\mathbb{Z}%
_{+})$ is defined by 
\begin{equation}
S(\alpha )(n):=\alpha (n+1)\;(\alpha \in \ell ^{\infty }(\mathbb{Z}%
_{+}),n\geq 0),  \label{sa}
\end{equation}%
then 
\begin{equation}
d\xi _{S(\alpha )}(s)=\frac{s}{\alpha _{0}^{2}}d\xi (s).  \label{sam}
\end{equation}

We now formally define the class $\mathcal{TC}$. \ First, we need some
notation: $\mathcal{M}_{1}:=\vee \{e_{k_{1},k_{2}}:k_{2}\geq 1\}$ and $%
\mathcal{N}_{1}:=\vee \{e_{k_{1},k_{2}}:k_{1}\geq 1\}$.

\begin{definition}
\label{tc} (i) \ The \textit{core} of a $2$-variable weighted shift $\mathbf{%
T}$ is $c(\mathbf{T}):=\mathbf{T}|_{\mathcal{M}_{1}\cap \mathcal{N}_{1}}$; 
\newline
(ii) $\mathbf{T}$ is said to be of \textit{tensor form} if $\mathbf{T}\cong
(I\otimes W_{\alpha },W_{\beta }\otimes I)$. $\ $(When $\mathbf{T}$ is
subnormal, this is equivalent to requiring that the Berger measure be a
Cartesian product $\xi \times \eta $); \newline
(iii) $\ \mathcal{TC}:=\{\mathbf{T}\in \mathfrak{H}_{0}:c(\mathbf{T})\text{
is of tensor form}\}$.
\end{definition}

\section{Main Results}

We now consider $2$-variable weighted shifts such as the one given by Figure 1(ii), where $W_{x}\equiv shift(x_{0},x_{1},\cdots )$
is subnormal with Berger measure $\xi _{x}$ and $W_{y}\equiv
shift(y_{0},y_{1},\cdots )$ is subnormal with Berger measure $\eta _{y}$. $\ 
$Further, let $W_{\alpha }\equiv shift(\alpha _{1},\alpha _{2},\cdots )$
(resp. $W_{\beta }\equiv shift(\beta _{1},\beta _{2},\cdots )$) be subnormal
with Berger measure $\xi $ (resp. $\eta $). \ By (\ref{sam}), and without
loss of generality, we will always assume that $\frac{1}{s}\in L^{1}(\xi )$
and $\frac{1}{t}\in L^{1}(\eta )$. \ We recall several notions introduced in %
\cite{CuYo1} and \cite{power}: (i) given a probability measure $\mu $ on $%
X\times Y\equiv \mathbb{R}_{+}\times \mathbb{R}_{+}$, with $\frac{1}{t}\in
L^{1}(\mu )$, the \textit{extremal measure} $\mu _{ext}$ (which is also a
probability measure) on $X\times Y$ is given by $d\mu _{ext}(s,t):=\frac{1}{%
t\left\| \frac{1}{t}\right\| _{L^{1}(\mu )}}d\mu (s,t)$; and (ii) given a
measure $\mu $ on $X\times Y$, the \textit{marginal measure} $\mu ^{X}$ $($%
resp. $\mu ^{Y})$ is given by $\mu ^{X}:=\mu \circ \pi _{X}^{-1}$ $($resp. $%
\mu ^{Y}:=\mu \circ \pi _{Y}^{-1})$, where $\pi _{X}:X\times Y\rightarrow X$ 
$($resp. $\pi _{Y}:X\times Y\rightarrow Y)$ is the canonical projection onto 
$X$ $($resp. $Y)$. \ Thus, $\mu ^{X}(E)=\mu (E\times Y)$, for every $%
E\subseteq X$ $($resp. $\mu ^{Y}(F)=\mu (X\times F)$, $F\subseteq Y)$. \
Observe that if $\mu $ is a probability measure, then so are $\mu ^{X}$ and $%
\mu ^{Y}$. \ 

For a measure $\mu $ with $\frac{1}{s}\in L^{1}(\mu )$, we write $d%
\widetilde{\mu }(s):=\frac{1}{s\left\| \frac{1}{s}\right\| _{L^{1}(\mu )}}%
d\mu (s)$. \ For example, 
\begin{equation}
d(\xi \times \eta )_{ext}(s,t)=\frac{1}{t\left\| \frac{1}{t}\right\|
_{L^{1}(\eta )}}d\xi (s)d\eta (t)=d\xi (s)d\widetilde{\eta }(t)  \label{eq1}
\end{equation}%
and $(\xi \times \eta )^{X}=\xi $. \ Finally, for an arbitrary $2$-variable
weighted shift $\mathbf{T}$, we shall let $\mathcal{R}_{ij}(\mathbf{T})$
denote the restriction of $\mathbf{T}$ to $\mathcal{M}_{i}\cap \mathcal{N}%
_{j},$ where $\mathcal{M}_{i}$(resp. $\mathcal{N}_{j})$ is the subspace of $%
\ell ^{2}(\mathbb{Z}_{+}^{2})$ spanned by the canonical orthonormal basis
associated with indices $\mathbf{k}=(k_{1},k_{2})$, where $k_{1}\geq 0$ and $%
k_{2}\geq i$ $(k_{1}\geq j$ and $k_{2}\geq 0)$, respectively. \ In
particular, we simply denote $\mathcal{M}\equiv \mathcal{M}_{1}$ and $%
\mathcal{N}\equiv \mathcal{N}_{1}$. \ It follows that $\mathcal{R}_{11}(%
\mathbf{T})=c(\mathbf{T})$, the core of $\mathbf{T}$. \ Assume that $c(%
\mathbf{T})$ is subnormal, with Berger measure $\xi \times \eta $. \ We let 
\begin{equation}
\psi :=\left( \eta _{y}\right) _{1}-a^{2}\left\| \frac{1}{s}\right\|
_{L^{1}(\xi )}\eta  \label{eqpsi}
\end{equation}%
and 
\begin{equation}
\varphi :=\xi _{x}-y_{0}^{2}\left\| \frac{1}{t}\right\| _{L^{1}(\psi
)}\delta _{0}-a^{2}y_{0}^{2}\left\| \frac{1}{s}\right\| _{L^{1}(\xi
)}\left\| \frac{1}{t}\right\| _{L^{1}(\eta )}\widetilde{\xi },  \label{eqphi}
\end{equation}%
where $\left( \eta _{y}\right) _{1}$ is the Berger measure of the subnormal
shift $shift(y_{1},y_{2},\cdots )$. \ Trivially, $\psi $ and $\varphi $ are
measures, but they may or may not be \textit{positive} measures. \ The
following result is a very special case of the Reconstruction-of-the-measure
Problem.

\begin{lemma}
\label{backext}(Subnormal backward extension of a $2$-variable weighted
shift \cite{CuYo1}) \ Consider the $2$-variable weighted shift whose weight
diagram is given in Figure 1(i). \ Assume that $%
\mathcal{R}_{10}(\mathbf{T})\equiv \mathbf{T}|_{\mathcal{M}}$ is subnormal,
with associated measure $\mu _{\mathcal{M}}$, and that $W_{0}\equiv
shift(\alpha _{00},\alpha _{10},\cdots )$ is subnormal with associated
measure $\nu $. \ Then $\mathbf{T}$ is subnormal if and only if\newline
(i) $\frac{1}{t}\in L^{1}(\mu _{\mathcal{M}})$;\newline
(ii) $\beta _{00}^{2}\leq (\left\| \frac{1}{t}\right\| _{L^{1}(\mu _{%
\mathcal{M}})})^{-1}$;\newline
(iii) $\beta _{00}^{2}\left\| \frac{1}{t}\right\| _{L^{1}(\mu _{\mathcal{M}%
})}(\mu _{\mathcal{M}})_{ext}^{X}\leq \nu $.

Moreover, if $\beta _{00}^{2}\left\| \frac{1}{t}\right\| _{L^{1}(\mu _{%
\mathcal{M}})}=1$ then $(\mu _{\mathcal{M}})_{ext}^{X}=\nu $. \ In the case
when $\mathbf{T}$ is subnormal, the Berger measure $\mu $ of $\mathbf{T}$ is
given by 
\begin{equation*}
d\mu (s,t)=\beta _{00}^{2}\left\| \frac{1}{t}\right\| _{L^{1}(\mu _{\mathcal{%
M}})}d(\mu _{\mathcal{M}})_{ext}(s,t)+(d\nu (s)-\beta _{00}^{2}\left\| \frac{%
1}{t}\right\| _{L^{1}(\mu _{\mathcal{M}})}d(\mu _{\mathcal{M}%
})_{ext}^{X}(s))d\delta _{0}(t).
\end{equation*}

\begin{figure}[th]
\includegraphics{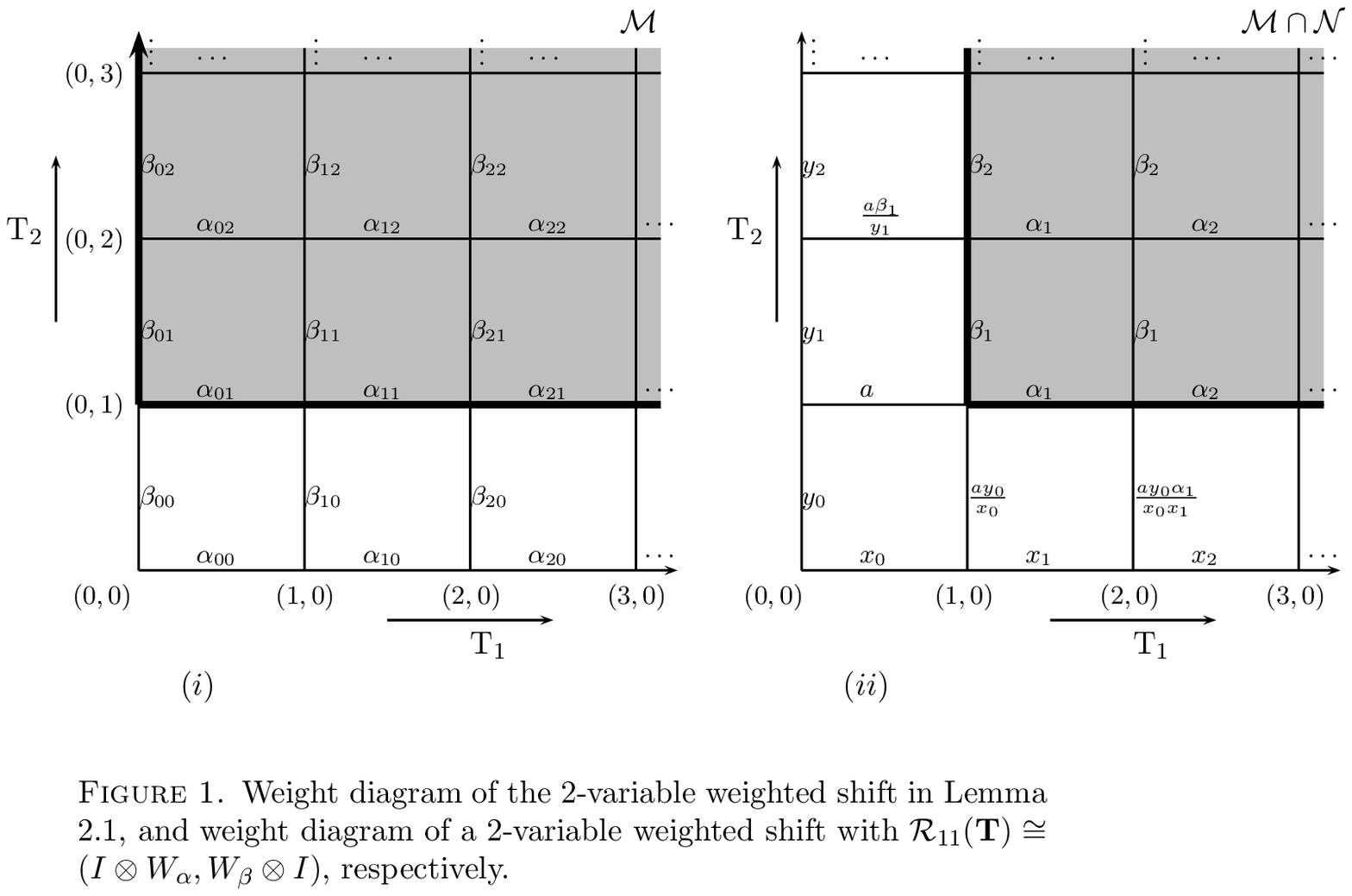}
\end{figure}

\end{lemma}

\begin{proposition}
\label{pro1}Let $\mathbf{T}\equiv \mathbf{(}T_{1},T_{2})\in \mathfrak{H}_{0}$
be the $2$-variable weighted shift whose weight diagram is given in 
Figure 1(ii). \ Then $\mathbf{T}|_{\mathcal{M}}\in 
\mathfrak{H}_{\infty }$ if and only if $\psi $ is a positive measure. \ In
this case, the Berger measure of $\mathbf{T}|_{\mathcal{M}}$ is 
\begin{equation*}
\mu _{\mathcal{M}}=a^{2}\left\| \frac{1}{s}\right\| _{L^{1}(\xi )}\widetilde{%
\xi }\times \eta +\delta _{0}\times \psi \text{.}
\end{equation*}
\end{proposition}

\begin{proof}
$(\Rightarrow )$ If $\mathbf{T}|_{\mathcal{M}}\in \mathfrak{H}_{\infty }$
then $\mathbf{T}|_{\mathcal{M\cap N}}\in \mathfrak{H}_{\infty }$ with Berger
measure $\mu _{\mathcal{M\cap N}}=\xi \times \eta $. \ Note that $\left\| 
\frac{1}{s}\right\| _{L^{1}(\mu _{\mathcal{M\cap N}})}=\left\| \frac{1}{s}%
\right\| _{L^{1}(\xi )}$. \ By Lemma \ref{backext}(iii), if we think of $%
\mathbf{T}|_{\mathcal{M}}$ as the backward extension of $\mathbf{T}|_{%
\mathcal{M\cap N}}$ (in the $s$ direction), we then have%
\begin{equation*}
\begin{tabular}{l}
$\left( \eta _{y}\right) _{1}\geq a^{2}\left\| \frac{1}{s}\right\|
_{L^{1}(\mu _{\mathcal{M\cap N}})}(\mu _{\mathcal{M\cap N}})_{ext}^{Y}$ \\ 
\\ 
$\Leftrightarrow \left( \eta _{y}\right) _{1}-a^{2}\left\| \frac{1}{s}%
\right\| _{L^{1}(\xi )}\eta \geq 0$ \\ 
\\ 
$\Leftrightarrow \psi \geq 0$.%
\end{tabular}%
\end{equation*}%
Thus, $\psi $ is a positive measure.\newline
$(\Leftarrow )$ If $\psi $ is a positive measure then $\phi :=a^{2}\left\| 
\frac{1}{s}\right\| _{L^{1}(\xi )}\widetilde{\xi }\times \eta +\delta
_{0}\times \psi $ is a well defined and positive measure. \ By a direct
calculation, we can see that 
\begin{equation*}
\left\{ 
\begin{tabular}{ll}
$\iint d\phi (s,t)=1,$ & $\text{if }k_{1}=0\text{ and }k_{2}=0$ \\ 
$\iint t^{k_{2}}d\phi (s,t)=y_{1}^{2}\cdots y_{k_{2}}^{2},$ & $\text{if }%
k_{1}=0\text{ and }k_{2}\geq 1$ \\ 
$\iint s^{k_{1}}d\phi (s,t)=a^{2}\alpha _{1}^{2}\cdots \alpha
_{k_{1}-1}^{2}, $ & $\text{if }k_{1}\geq 1\text{ and }k_{2}=0$ \\ 
$\iint s^{k_{1}}t^{k_{2}}d\phi (s,t)=a^{2}\alpha _{1}^{2}\cdots \alpha
_{k_{1}-1}^{2}\beta _{1}^{2}\cdots \beta _{k_{2}}^{2},$ & $\text{if }%
k_{1}\geq 1\text{ and }k_{2}\geq 1$%
\end{tabular}%
\right\} =\gamma _{\mathbf{k}}(\mathbf{T}|_{\mathcal{M}}),
\end{equation*}%
where, for notational convenience, we set $\alpha _{0}:=1$. \ Therefore, $%
\phi $ interpolates all moments of $\mathbf{T}|_{\mathcal{M}}$, so $\mathbf{T%
}|_{\mathcal{M}}\in \mathfrak{H}_{\infty }$ and $\mu _{\mathcal{M}}=\phi
\equiv a^{2}\left\| \frac{1}{s}\right\| _{L^{1}(\xi )}\widetilde{\xi }\times
\eta +\delta _{0}\times \psi $.
\end{proof}

We now have:

\begin{theorem}
\label{thm1}Let $\mathbf{T}\equiv \mathbf{(}T_{1},T_{2})\in \mathfrak{H}_{0}$
be the $2$-variable weighted shift whose weight diagram is given in 
Figure 1(ii). \ Then $\mathbf{T}\in \mathfrak{H}_{\infty }$
if and only if $\psi $ and $\varphi $ are positive measures.
\end{theorem}

\begin{proof}
$(\Leftarrow )$ It suffices to find a probability measure $\mu $ satisfying%
\begin{equation*}
\gamma _{\mathbf{k}}(\mathbf{T})=\int s^{k_{1}}t^{k_{2}}d\mu (s,t)\;\;(\text{%
all }\mathbf{k}\equiv (k_{1},k_{2})\in \mathbb{Z}_{+}^{2}).
\end{equation*}%
\ Let 
\begin{equation*}
\mu :=\varphi \times \left( \delta _{0}-\widetilde{\eta }\right)
+y_{0}^{2}\left\| \frac{1}{t}\right\| _{L^{1}(\psi )}\delta _{0}\times
\left( \widetilde{\psi }-\widetilde{\eta }\right) +\xi _{x}\times \widetilde{%
\eta }.
\end{equation*}%
Clearly, $\mu $ is well defined. \ Observe that%
\begin{eqnarray}
\mu  &=&\varphi \times \left( \delta _{0}-\widetilde{\eta }\right)
+y_{0}^{2}\left\| \frac{1}{t}\right\| _{L^{1}(\psi )}\delta _{0}\times
\left( \widetilde{\psi }-\widetilde{\eta }\right) +\xi _{x}\times \widetilde{%
\eta }  \notag \\
&&  \notag \\
&=&(\xi _{x}-\varphi -y_{0}^{2}\left\| \frac{1}{t}\right\| _{L^{1}(\psi
)}\delta _{0})\times \widetilde{\eta }+y_{0}^{2}\left\| \frac{1}{t}\right\|
_{L^{1}(\psi )}\delta _{0}\times \widetilde{\psi }+\varphi \times \delta _{0}
\notag \\
&&  \notag \\
&=&a^{2}y_{0}^{2}\left\| \frac{1}{s}\right\| _{L^{1}(\xi )}\left\| \frac{1}{t%
}\right\| _{L^{1}(\eta )}\widetilde{\xi }\times \widetilde{\eta }%
+y_{0}^{2}\left\| \frac{1}{t}\right\| _{L^{1}(\psi )}\delta _{0}\times 
\widetilde{\psi }+\varphi \times \delta _{0}\text{.}  \label{eq2}
\end{eqnarray}%
Since we are assuming that $\psi $ and $\varphi $ are positive measures, it
follows from \ref{eq2} that $\mu $ is also positive. \ Furthermore, observe
that 
\begin{equation}
\left\| \frac{1}{t}\right\| _{L^{1}(\psi )}=\left\| \frac{1}{t}\right\|
_{L^{1}(\left( \eta _{y}\right) _{1})}-a^{2}\left\| \frac{1}{s}\right\|
_{L^{1}(\xi )}\left\| \frac{1}{t}\right\| _{L^{1}(\eta )}  \label{normofpsi}
\end{equation}%
and%
\begin{equation*}
\left\| \frac{1}{t}\right\| _{L^{1}(\psi )}\widetilde{\psi }=\left\| \frac{1%
}{t}\right\| _{L^{1}(\left( \eta _{y}\right) _{1})}\widetilde{\left( \eta
_{y}\right) }_{1}-a^{2}\left\| \frac{1}{s}\right\| _{L^{1}(\xi )}\left\| 
\frac{1}{t}\right\| _{L^{1}(\eta )}\widetilde{\eta }\text{.}
\end{equation*}%
Thus,%
\begin{eqnarray*}
\iint d\mu (s,t) &=&y_{0}^{2}\left\| \frac{1}{t}\right\| _{L^{1}(\psi )}\int
d\widetilde{\psi }(t)-y_{0}^{2}\left\| \frac{1}{t}\right\| _{L^{1}(\psi )}+1
\\
&& \\
&=&1\text{ \ (since }\widetilde{\psi }\text{ is a probability measure).}
\end{eqnarray*}%
Therefore, $\mu $ is a probability measure.

Next, observe that%
\begin{equation}
\varphi ([0,+\infty ))=\int d\varphi (s)=1-y_{0}^{2}\left\| \frac{1}{t}%
\right\| _{L^{1}(\psi )}-a^{2}y_{0}^{2}\left\| \frac{1}{s}\right\|
_{L^{1}(\xi )}\left\| \frac{1}{t}\right\| _{L^{1}(\eta )}  \label{eq20}
\end{equation}%
and, for $k_{2}\geq 1$,%
\begin{equation}
\left\| \frac{1}{t}\right\| _{L^{1}(\psi )}\int t^{k_{2}}d\widetilde{\psi }%
(t)=\int t^{k_{2}-1}d\left( \eta _{y}\right) _{1}(t)-a^{2}\left\| \frac{1}{s}%
\right\| _{L^{1}(\xi )}\int t^{k_{2}-1}d\eta (t)\text{.}  \label{eq21}
\end{equation}%
Thus, if $k_{1}=0$ and $k_{2}\geq 1$, we have 
\begin{eqnarray*}
\iint t^{k_{2}}d\mu (s,t) &=&-\int d\varphi (s)\int t^{k_{2}}d\widetilde{%
\eta }(t)+y_{0}^{2}\left\| \frac{1}{t}\right\| _{L^{1}(\psi )}\int
t^{k_{2}}d\left( \widetilde{\psi }-\widetilde{\eta }\right) (t)+\int
t^{k_{2}}d\widetilde{\eta }(t) \\
&& \\
&=&a^{2}y_{0}^{2}\left\| \frac{1}{s}\right\| _{L^{1}(\xi )}\left\| \frac{1}{t%
}\right\| _{L^{1}(\eta )}\int t^{k_{2}}d\widetilde{\eta }(t)+y_{0}^{2}\left%
\| \frac{1}{t}\right\| _{L^{1}(\psi )}\int t^{k_{2}}d\widetilde{\psi }(t) \\
&&\text{(by (\ref{eq20}))} \\
&=&y_{0}^{2}\int t^{k_{2}-1}d\left( \eta _{y}\right) _{1}(t)=y_{0}^{2}\cdots
y_{k_{2}-1}^{2} \\
&&\text{(by (\ref{eq21})).}
\end{eqnarray*}%
If $k_{1}\geq 1$ and $k_{2}=0$, we have 
\begin{eqnarray*}
\iint s^{k_{1}}d\mu (s,t) &=&\iint s^{k_{1}}d\varphi (s)d\left( \delta _{0}-%
\widetilde{\eta }\right) (t)+\iint s^{k_{1}}d\xi _{x}(s)d\widetilde{\eta }(t)
\\
&& \\
&=&\int s^{k_{1}}d\xi _{x}(s) \\
&& \\
&=&x_{0}^{2}\cdots x_{k_{1}-1}^{2}\text{.}
\end{eqnarray*}%
Finally, if $k_{1}\geq 1$ and $k_{2}\geq 1$, we have 
\begin{eqnarray*}
\iint s^{k_{1}}t^{k_{2}}d\mu (s,t) &=&\iint s^{k_{1}}t^{k_{2}}d\varphi
(s)d\left( \delta _{0}-\widetilde{\eta }\right) (t)+\iint
s^{k_{1}}t^{k_{2}}d\xi _{x}(s)d\widetilde{\eta }(t) \\
&& \\
&=&\iint s^{k_{1}}t^{k_{2}}(d\xi _{x}(s)-d\varphi (s))d\widetilde{\eta }(t)
\\
&& \\
&=&a^{2}y_{0}^{2}\left\| \frac{1}{s}\right\| _{L^{1}(\xi )}\left\| \frac{1}{t%
}\right\| _{L^{1}(\eta )}\iint s^{k_{1}}t^{k_{2}}d\widetilde{\xi }(s)d%
\widetilde{\eta }(t) \\
&& \\
&=&a^{2}y_{0}^{2}\alpha _{1}^{2}\cdots \alpha _{k_{1}-1}^{2}\beta
_{1}^{2}\cdots \beta _{k_{2}-1}^{2},
\end{eqnarray*}%
where, for notational convenience, we set $\alpha _{0}:=1$ and $\beta _{0}:=1
$. \ Thus,%
\begin{equation*}
\left\{ 
\begin{tabular}{ll}
$\iint d\mu (s,t)=1,$ & $\text{if }k_{1}=0\text{ and }k_{2}=0$ \\ 
$\iint t^{k_{2}}d\mu (s,t)=y_{0}^{2}\cdots y_{k_{2}-1}^{2},$ & $\text{if }%
k_{1}=0\text{ and }k_{2}\geq 1$ \\ 
$\iint s^{k_{1}}d\mu (s,t)=x_{0}^{2}\cdots x_{k_{1}-1}^{2},$ & $\text{if }%
k_{1}\geq 1\text{ and }k_{2}=0$ \\ 
$\iint s^{k_{1}}t^{k_{2}}d\mu (s,t)=a^{2}y_{0}^{2}\alpha _{1}^{2}\cdots
\alpha _{k_{1}-1}^{2}\beta _{1}^{2}\cdots \beta _{k_{2}-1}^{2},$ & $\text{if 
}k_{1}\geq 1\text{ and }k_{2}\geq 1$%
\end{tabular}%
\right\} =\gamma _{\mathbf{k}}(\mathbf{T})\text{.}
\end{equation*}%
Therefore, $\mu $ interpolates all moments of $\mathbf{T}$, so $\mathbf{T}$
must be subnormal, with Berger measure $\mu $.

$(\Rightarrow )$ Assume that $\mathbf{T}$ is subnormal with Berger measure $%
\mu $. \ Then $\mathbf{T}|_{\mathcal{M}}$ is also subnormal, and by
Proposition \ref{pro1} we can see that $\psi $ is a positive measure.\ \ We
then have 
\begin{eqnarray*}
\left\| \frac{1}{t}\right\| _{L^{1}(\mu _{\mathcal{M}})} &=&\int \frac{1}{t}%
d\mu _{\mathcal{M}}(s,t) \\
&& \\
&=&a^{2}\left\| \frac{1}{s}\right\| _{L^{1}(\xi )}\left\| \frac{1}{t}%
\right\| _{L^{1}(\eta )}+\int \frac{1}{t}d\left( \eta _{y}\right)
_{1}(t)-a^{2}\left\| \frac{1}{s}\right\| _{L^{1}(\xi )}\left\| \frac{1}{t}%
\right\| _{L^{1}(\eta )} \\
&& \\
&=&\left\| \frac{1}{t}\right\| _{L^{1}(\left( \eta _{y}\right) _{1})}.
\end{eqnarray*}%
Since $\mu _{\mathcal{M}}=a^{2}\left\| \frac{1}{s}\right\| _{L^{1}(\xi )}%
\widetilde{\xi }\times \eta +\delta _{0}\times \psi \text{,}$ we get 
\begin{eqnarray*}
\left\| \frac{1}{t}\right\| _{L^{1}(\left( \eta _{y}\right) _{1})}d(\mu _{%
\mathcal{M}})_{ext}(s,t) &=&\left\| \frac{1}{t}\right\| _{L^{1}(\mu _{%
\mathcal{M}})}d(\mu _{\mathcal{M}})_{ext}(s,t) \\
&& \\
&=&\left\| \frac{1}{t}\right\| _{L^{1}(\mu _{\mathcal{M}})}d\left\{
a^{2}\left\| \frac{1}{s}\right\| _{L^{1}(\xi )}\widetilde{\xi }\times \eta
+\delta _{0}\times \psi \right\} _{ext}(s,t) \\
&& \\
&=&a^{2}\left\| \frac{1}{s}\right\| _{L^{1}(\xi )}d\widetilde{\xi }(s)\frac{%
d\eta (t)}{t}+d\delta _{0}(s)\frac{d\psi (t)}{t} \\
&& \\
&=&a^{2}\left\| \frac{1}{s}\right\| _{L^{1}(\xi )}d\widetilde{\xi }(s)\frac{%
d\eta (t)}{t}+d\delta _{0}(s)[\frac{d\left( \eta _{y}\right) _{1}(t)}{t}%
-a^{2}\left\| \frac{1}{s}\right\| _{L^{1}(\xi )}\frac{d\eta (t)}{t}].
\end{eqnarray*}%
It follows that 
\begin{eqnarray*}
\left\| \frac{1}{t}\right\| _{L^{1}(\left( \eta _{y}\right) _{1})}(\mu _{%
\mathcal{M}})_{ext}^{X} &=&\left\| \frac{1}{t}\right\| _{L^{1}(\left( \eta
_{y}\right) _{1})}\int_{Y}d(\mu _{\mathcal{M}})_{ext}(\cdot ,t) \\
&& \\
&=&a^{2}\left\| \frac{1}{s}\right\| _{L^{1}(\xi )}\left\| \frac{1}{t}%
\right\| _{L^{1}(\eta )}\widetilde{\xi }+\left( \left\| \frac{1}{t}\right\|
_{L^{1}(\left( \eta _{y}\right) _{1})}-a^{2}\left\| \frac{1}{s}\right\|
_{L^{1}(\xi )}\left\| \frac{1}{t}\right\| _{L^{1}(\eta )}\right) \delta _{0}.
\end{eqnarray*}%
Now recall that $\beta _{00}^{2}=y_{0}^{2}$, so from Lemma \ref{backext}%
(iii) we obtain%
\begin{equation*}
\xi _{x}\geq a^{2}y_{0}^{2}\left\| \frac{1}{s}\right\| _{L^{1}(\xi )}\left\| 
\frac{1}{t}\right\| _{L^{1}(\eta )}\widetilde{\xi }+y_{0}^{2}\left\| \frac{1%
}{t}\right\| _{L^{1}(\psi )}\delta _{0}\;\;\text{(using (\ref{normofpsi})}.
\end{equation*}%
Thus, $\varphi $ is a positive measure, as desired.
\end{proof}

\begin{remark}
The proof of Theorem \ref{thm1} gives a concrete formula for the Berger
measure of $\mathbf{T}$, namely 
\begin{equation*}
\mu =\varphi \times \left( \delta _{0}-\widetilde{\eta }\right)
+y_{0}^{2}\left\| \frac{1}{t}\right\| _{L^{1}(\psi )}\delta _{0}\times
\left( \widetilde{\psi }-\widetilde{\eta }\right) +\xi _{x}\times \widetilde{%
\eta }.
\end{equation*}
\end{remark}

In Proposition \ref{pro1} and Theorem \ref{thm1} we noted that $\psi $
(resp. $\varphi $) is a linear combination of $\left( \eta _{y}\right) _{1}$
and $\eta $ (resp. $\xi _{x}$, $\delta _{0}$ and $\widetilde{\xi }$), where $%
\left( \eta _{y}\right) _{1}$ and $\eta $ $($resp. $\xi _{x}$, $\delta _{0}$
and $\widetilde{\xi })$ are the Berger measures of the subnormal $1$%
-variable weighted shifts in the vertical (resp. horizontal) slices of $%
\mathbf{T}$. \ Thus, the following conjecture for $2$-variable weighted
shifts seems natural.

\begin{conjecture}
Let $\mathbf{T}\equiv \mathbf{(}T_{1},T_{2})\in \mathfrak{H}_{0}$ be the $2$%
-variable weighted shift whose weight diagram is given by 
Figure 1(i). \ Then the subnormality of $\mathbf{T}$ is
determined by a countable collection of inequalities $\{\psi _{k}\geq 0\}$,
where each measure $\psi _{k}$ is a linear combination of Berger measures
associated to the $1$-variable weighted shifts in vertical or horizontal
slices of $\mathbf{T}$.
\end{conjecture}

\section{Application to Flat $2$-variable Weighted Shifts}

We can now give a concrete condition for the subnormality of flat $2$%
-variable weighted shifts $\mathbf{T\equiv (}T_{1},T_{2}).$ \ Recall that $%
\mathbf{T}\equiv (T_{1},T_{2})$ is called \textit{horizontally flat} if $%
\alpha _{(k_{1},k_{2})}=\alpha _{(1,1)}$ for all $k_{1},k_{2}\geq 1$, and 
\textit{vertically flat} if $\beta _{(k_{1},k_{2})}=\beta _{(1,1)}$ for all $%
k_{1},k_{2}\geq 1$ \cite{propagation}. \ If $\mathbf{T}$ is horizontally and
vertically flat, then $\mathbf{T}$ is simply called \textit{flat} (see
Figure 2). \ It is straightforward to prove that $\mathbf{T}$ is
flat if and only if $\mathbf{T}\in \mathcal{TC}$, with $\xi $ and $\eta $ $1$%
-atomic. \ Without loss of generality, we can always assume that $\xi
=\delta _{1}$. \ 

\begin{figure}[th]
\includegraphics{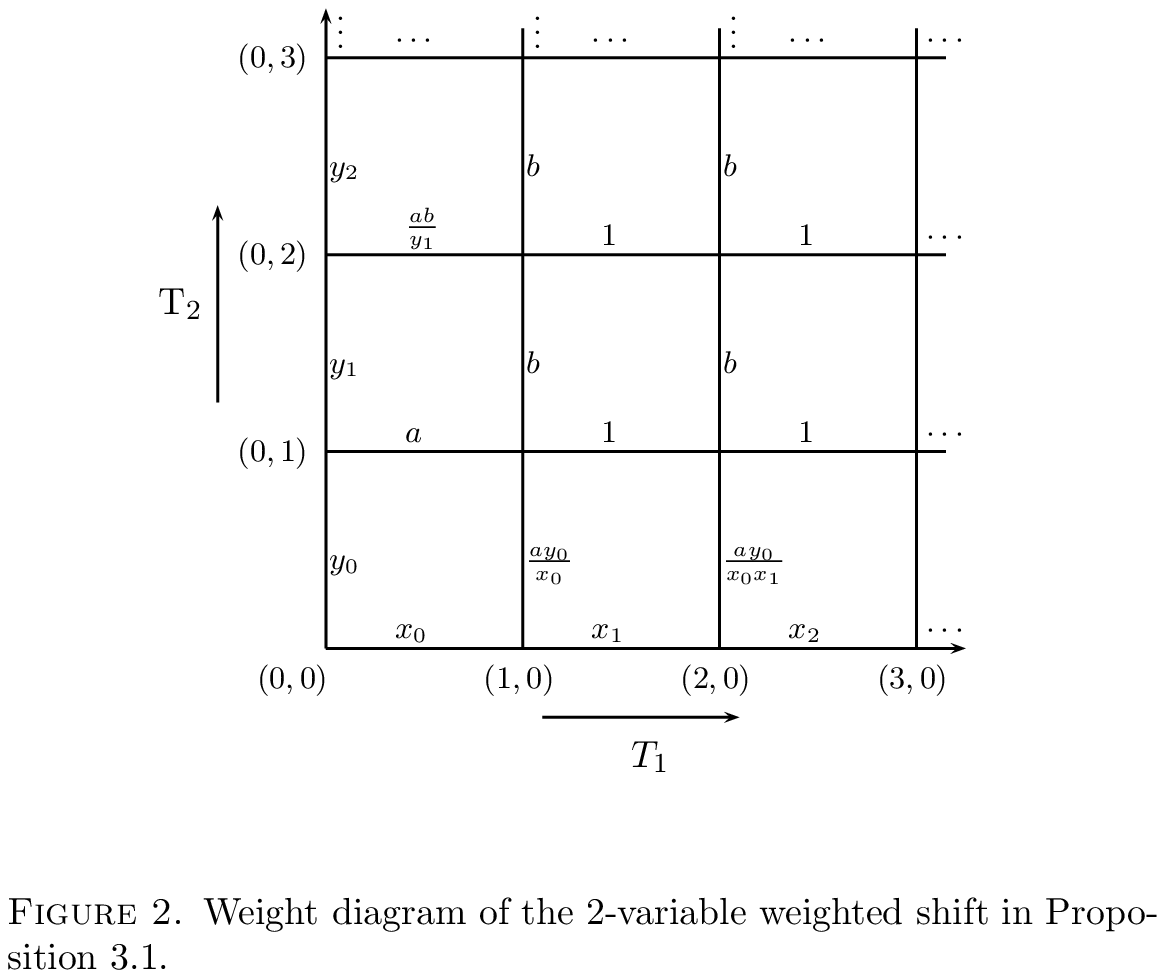}
\end{figure}

Now recall that for $0<\alpha <\beta $, $shift(\alpha ,\beta ,\beta ,\cdots
) $ is subnormal with Berger measure $(1-\frac{\alpha ^{2}}{\beta ^{2}}%
)\delta _{0}+\frac{\alpha ^{2}}{\beta ^{2}}\delta _{\beta ^{2}}$. \ Finally,
we know from \cite[Theorem 3.3]{CuYo2} and \cite[Section 5]{propagation}
that if $\mathbf{T}$ is flat and subnormal then $\xi _{x}$ and $\eta _{y}$
have the form 
\begin{equation}
\begin{tabular}{l}
$\xi _{x}=p\delta _{0}+q\delta _{1}+[1-(p+q)]\rho $ \\ 
$\eta _{y}=\ell \delta _{0}+m\delta _{b^{2}}+[1-(\ell +m)]\sigma $,%
\end{tabular}
\label{eta}
\end{equation}%
where $0<p,q,\ell ,m<1,$ $p+q\leq 1$, $\ell +m\leq 1$, and $\rho $ and $%
\sigma $ are probability measures with $\rho (\{0\}\cup \{1\})=0$, $\sigma
(\{0\}\cup \{b^{2}\})=0$.

We are now ready to present

\begin{proposition}
\label{ex1} Consider the $2$-variable weighted shift $\mathbf{T}\in 
\mathfrak{H}_{0}$ whose weight diagram is given in Figure 2. \ The
following statements are equivalent.\newline
(i) $\ \mathbf{T}\in \mathfrak{H}_{\infty }$;\newline
(ii) \ $\psi $ and $\varphi $ are positive measures;\newline
(iii) \ $\frac{b}{a}\sqrt{m}\geq y_{0}$ and 
\begin{equation}
\xi _{x}\geq y_{0}^{2}\left\{ \left( \left\| \frac{1}{t}\right\|
_{L^{1}(\left( \eta _{y}\right) _{1})}-\frac{a^{2}}{b^{2}}\right) \delta
_{0}+\frac{a^{2}}{b^{2}}\delta _{1}\right\} .  \label{eq30}
\end{equation}%
\newline
Moreover, when $\mathbf{T}$ is subnormal, its Berger measure is given as 
\begin{equation}
\mu =\varphi \times \left( \delta _{0}-\delta _{b^{2}}\right)
+y_{0}^{2}\left\| \frac{1}{t}\right\| _{L^{1}(\psi )}\delta _{0}\times
\left( \widetilde{\psi }-\delta _{b^{2}}\right) +\delta _{1}\times \delta
_{b^{2}}  \label{eq4}
\end{equation}
\end{proposition}

\begin{proof}
(i) $\Rightarrow $ (ii): This is straightforward from Theorem \ref{thm1}.%
\newline
(ii) $\Rightarrow $ (iii): In (\ref{eqpsi}) and (\ref{eqphi}), observe that $%
\widetilde{\xi }=\delta _{1}$, $\widetilde{\eta }=\delta _{b^{2}}$, $\left\| 
\frac{1}{s}\right\| _{L^{1}(\xi )}=1$, $\left\| \frac{1}{t}\right\|
_{L^{1}(\eta )}=\frac{1}{b^{2}}$ and $\left\| \frac{1}{t}\right\|
_{L^{1}(\psi )}=\left\| \frac{1}{t}\right\| _{L^{1}(\left( \eta _{y}\right)
_{1})}-\frac{a^{2}}{b^{2}}$. $\ $Thus, we have 
\begin{eqnarray}
\psi &=&\left( \eta _{y}\right) _{1}-a^{2}\left\| \frac{1}{s}\right\|
_{L^{1}(\xi )}\eta  \notag \\
&&  \notag \\
&=&\frac{mb^{2}}{y_{0}^{2}}\delta _{b^{2}}+\frac{[1-(\ell +m)]}{y_{0}^{2}}%
\sigma _{1}-a^{2}\delta _{b^{2}}  \label{eq5}
\end{eqnarray}%
(where $d\sigma _{1}(t):=td\sigma (t)$) and%
\begin{eqnarray}
\varphi &=&\xi _{x}-y_{0}^{2}\left\| \frac{1}{t}\right\| _{L^{1}(\psi
)}\delta _{0}-a^{2}y_{0}^{2}\left\| \frac{1}{s}\right\| _{L^{1}(\xi
)}\left\| \frac{1}{t}\right\| _{L^{1}(\eta )}\widetilde{\xi }  \notag \\
&&  \notag \\
&=&\xi _{x}-y_{0}^{2}\left\{ \left( \left\| \frac{1}{t}\right\|
_{L^{1}(\left( \eta _{y}\right) _{1})}-\frac{a^{2}}{b^{2}}\right) \delta
_{0}+\frac{a^{2}}{b^{2}}\delta _{1}\right\} \text{.}  \label{eq6}
\end{eqnarray}%
Since we are assuming that $\psi $ and $\varphi $ are positive measures, it
follows from (\ref{eq5}) and (\ref{eq6}) that 
\begin{equation*}
\frac{b}{a}\sqrt{m}\geq y_{0}\text{ and }\xi _{x}\geq y_{0}^{2}\left\{
\left( \left\| \frac{1}{t}\right\| _{L^{1}(\left( \eta _{y}\right) _{1})}-%
\frac{a^{2}}{b^{2}}\right) \delta _{0}+\frac{a^{2}}{b^{2}}\delta
_{1}\right\} ,
\end{equation*}%
as desired.\newline
(iii) $\Rightarrow $ (i): \ Let $\omega :=a^{2}\delta _{1}\times \delta
_{b^{2}}+\delta _{0}\times \psi $. \ Then $\omega $ is well defined, and by
the formula for $\psi $ (given in (\ref{eq5})) and the condition $\frac{b}{x}%
\sqrt{m}\geq y_{0}$, we see at once that $\omega \geq 0$. \ Furthermore, $%
\omega $ is the Berger measure of $\mathbf{T}|_{\mathcal{M}}$, so that $%
\mathbf{T}|_{\mathcal{M}}$ is subnormal. \ We now wish to apply Lemma \ref%
{backext}. \ Note that%
\begin{equation*}
d\omega _{ext}(s,t)=\frac{1}{t\left\| \frac{1}{t}\right\| _{L^{1}(\omega )}}%
\left\{ a^{2}d\delta _{1}(s)d\delta _{b^{2}}(t)+d\delta _{0}(s)d\psi
(t)\right\} ,
\end{equation*}%
so that%
\begin{equation*}
\omega _{ext}^{X}=\frac{1}{\left\| \frac{1}{t}\right\| _{L^{1}(\eta _{1})}}%
\left\{ \left( \left\| \frac{1}{t}\right\| _{L^{1}(\eta _{1})}-\frac{a^{2}}{%
b^{2}}\right) \delta _{0}+\frac{a^{2}}{b^{2}}\delta _{1}\right\} .
\end{equation*}%
We now see that the conditions (i), (ii) and (iii) in Lemma \ref{backext}
are satisfied, and therefore $\mathbf{T\in }\mathfrak{H}_{\infty }$. $\ $

Finally, since in this case we have $\widetilde{\xi }=\delta _{1}$ and $%
\widetilde{\eta }=\delta _{b^{2}}$, Theorem \ref{thm1} readily implies (\ref%
{eq4}). \ 
\end{proof}

\end{document}